\documentclass[12pt]{amsart}
 
\usepackage{tikz}
\usetikzlibrary{cd}

\usepackage{amsfonts,amssymb}

\usepackage{pinlabel,url}
\usepackage[all]{xy}
\usepackage{tikz}
\usepackage{graphicx}
\usepackage{xcolor}
\usepackage{hyperref}
\usepackage{comment}

\usepackage[nameinlink]{cleveref}

\newtheorem{theorem}{Theorem}[section]

\theoremstyle{definition}

\newtheorem*{remark}{Remark}

\newtheorem{question}[theorem]{Question}

\numberwithin{equation}{section}

\newcommand{\Mod}{{\rm Mod}}
\newcommand{\Out}{{\rm Out}}
\newcommand{\QI}{{\rm QI}}

\newcounter{mcomments}

\newcounter{tcomments}

\title[$\mathcal A(S_{g,1})$ is not quasi-isometric to $\mathcal{S}_{2g}$]{The arc complex is not quasi-isometric \\ to the sphere complex}

\author[T. Aougab and M. Loving]{Tarik Aougab and Marissa Loving}

\date{\today}

\begin{document}

\begin{abstract}
  We show that the arc complex $\mathcal{A}(S_{g,1})$ is not quasi-isometric to the sphere complex $\mathcal S_{2g}$ associated to the double of a genus $2g$ handlebody. Along the way, we present a simple proof that $\mathcal{A}(S_{g,1})$ is quasi-isometrically rigid. 
\end{abstract}

\maketitle

%%%%%%% Add notes here %%%%%%%

\section{Introduction}
\label{sec:intro}

Let $S = S_{g,p}$ be a connected, orientable, finite-type surface of genus $g$ with $p$ punctures. The \textit{arc complex} of $S$, denoted $\mathcal A(S)$, is the simplicial complex whose $k$-simplices correspond to collections of $k+1$ isotopy classes of properly embedded essential arcs that can be realized simultaneously disjointly on $S$. 

Let $H_{g}$ denote a handlebody of genus $g$, which is topologically the connect sum of $g$ copies of a solid torus. Then the \textit{double} of $H_{g}$, denoted $D H_g$, is the closed $3$-manifold obtained by gluing two copies of $H_{g}$ over their boundary by the identity map; this is topologically the connect sum of $g$ copies of $S^{2} \times S^{1}$. The \textit{sphere complex of rank $g$}, $\mathcal{S}_{g}$, is the simplicial complex whose $k$-simplices correspond to collections of $k+1$ isotopy classes of embedded essential spheres in the double that can be realized simultaneously disjointly.

There exists a natural map 
\[ \phi: \mathcal{A}(S_{g,1}) \rightarrow \mathcal{S}_{2g}\]
defined as follows. Consider the product of $S_{g,1}$ (where we blow up the puncture to a boundary component) with $[0,1]$; the result up to homeomorphism is a genus $2g$ handlebody $H_{2g}$, and each properly embedded arc in $S_{g,1}$ becomes a meridian. Doubling $S_{g,1} \times [0,1] \cong H_{2g}$ therefore associates each arc on $S_{g,1}$ to an essential sphere. Hamendst\"adt--Hensel \cite{HamenstadtHensel-2014} show that $\phi$ is a quasi-isometric embedding. Thus, it is natural to ask the following:

\begin{question} \label{question:qi}
    Is $\mathcal{A}(S_{g,1})$ (abstractly) quasi-isometric to $\mathcal{S}_{2g}$?
\end{question}

The difficulty in answering this question stems from the lack of known invariants of $\mathcal{S}_{2g}$ that would help in deciding its quasi-isometry type. The purpose of this note is to provide a negative answer to \Cref{question:qi}. 

\begin{theorem} \label{thm:main} The complexes $\mathcal A(S_{g,1})$ and $\mathcal{S}_{2g}$ are not quasi-isometric. 
\end{theorem}

We think of this result as being of a similar flavor to recent work of Aramayona-Parlier-Webb \cite{APW25} which provides a criterion for establishing that many of the standard combinatorial complexes associated to a surface or to a free group (including the arc complex and the sphere complex) are not quasi-isometric to the curve complex of a non-sporadic surface. Coupling our result with theirs, one concludes that $\mathcal{A}(S_{g,1}), \mathcal{S}_{2g}$, and any curve complex $\mathcal{C}(Y)$ are pairwise non quasi-isometric. 

\subsection{Background and strategy.}

The arc complex is an essential combinatorial tool for studying the mapping class group $\Mod(S)$ which acts on $\mathcal{A}(S)$ simplicially. It was shown by Masur--Schleimer that $\mathcal{A}(S)$ is Gromov hyperbolic \cite{MasurSchleimer2013} and its geometry connects deeply to the study of $3$-manifolds that fiber over $S^{1}$ \cite{Strenner2023}. Furthermore, Disarlo proved that $\mathcal A(S)$ satisfies an Ivanov-type rigidity; the full group of simplicial automorphism of $\mathcal{A}(S)$ is canonically isomorphic to the extended mapping class group $\Mod^{\pm}(S)$ \cite[Theorem~1.2]{Disarlo2015}.

A homeomorphism of $D H_g$ sends essential spheres to essential spheres and preserves disjointness; thus it induces a map of $\mathcal{S}_{g}$ depending only on the isotopy class of the homeomorphism. It follows that $\mathcal{S}_{g}$ supports an action of $\mbox{Out}(F_{2g})$ by simplicial automorphisms. In fact, Aramayona-Souto show that the full group of simplicial automorphisms of $\mathcal S_g$ is canonically isomorphic to $\mbox{Out}(F_{2g})$ for $g \geq 3$ \cite{AramayonaSouto2011}. 

We prove the main theorem by leveraging the difference in size between the groups of quasi-isometries of the two complexes, and appealing to the fact that the mapping class group is co-Hopfian. A crucial ingredient in the proof is the following description of the group of quasi-isometries for $\mathcal{A}(S_{g,1})$:

\begin{theorem} \label{thm:rigidity}
Every quasi-isometry of $\mathcal A(S_{g,1})$ is bounded distance from a simplicial automorphism.
\end{theorem}

\section{Proof of the Main Theorem}

Note that \Cref{thm:rigidity} can be quantified as follows. 

\begin{theorem} \label{thm:quantified} There exists $D$ so that if 
\[ \eta: \mathcal{A}(S) \rightarrow \mathcal{A}(S)\]
is a quasi-isometry, then there exists a simplicial automorphism $f$ of $\mathcal{A}(S)$ so that for any arc $\alpha$, 
\[ d_{\mathcal{A}}(f(\alpha), \eta(\alpha)) \leq D \]    
\end{theorem}

This motivates a relation on the set of quasi-isometries of a space $X$. Given $f,g:X \rightarrow X$ quasi-isometries, we write $f \sim g$ when they are bounded distance apart, namely, there is some finite $D>0$ so that the images of any point under $f$ and $g$ have distance at most $D$ from one another. Then $\sim$ is an equivalence relation, and the set of equivalence classes admits a group structure with respect to composition-- this is the group of quasi-isometries of $X$, denoted $QI(X)$. 

An important consequence of Theorem \ref{thm:rigidity}, together with \cite[Theorem~1.2]{Disarlo2015} is the fact that $QI(\mathcal{A}(S))$ is isomorphic to $\Mod^{\pm}(S)$, which we will make use of in our proof of \Cref{thm:main}. 

\vspace{1 mm}

\begin{proof}[Proof of \Cref{thm:main}] Assume there is a quasi-isometry 
\[ f: \mathcal{A}(S_{g,1}) \rightarrow \mathcal{S}_{2g}.\]
Then $f$ induces an isomorphism $f_{\ast}$ between the quasi-isometry groups of $\mathcal A(S_{g,1})$ and $\mathcal S_{2g}$, respectively. Note that $\Out(F_{2g})$ is contained as a subgroup of $\QI(\mathcal S_{2g})$, since, for any simplicial complex $X$, $\QI(X)$ always contains the group of simplicial automorphisms of $X$. Since $\pi_{1}(S_{g,1})$ is isomorphic to $F_{2g}$, then $\Mod^{\pm}(S_{g,1})$ is contained as a subgroup of $\Out(F_{2g})$. In particular we have $\Mod^{\pm}(S_{g,1}) < \Out(F_{2g}) < \QI(\mathcal S_{2g})$. Denote by $G$ this copy of $\Mod^{\pm}(S_{g,1})$ in $\QI(\mathcal S_{2g})$.

\medskip

\noindent {\bf Claim:} $G$ is a proper subgroup of $\QI( \mathcal S_{2g})$.

\begin{proof}[Proof of Claim] It suffices to show that there exists a quasi-isometry of the sphere complex that is not bounded distance from any mapping class. First note that no element of $G < \Out(F_{2g})$ is an atoroidal outer automorphism of $F_{2g}$ because it necessarily fixes the boundary component of $S_{g,1}$. Next consider a fully irreducible and atoroidal element $\varphi \in \Out(F_{2g})$. It necessarily acts with north-south dynamics on the compactified outer space $\overline{CV}_n$ \cite{LevittLustig}; associated to $\varphi$ are exactly two fixed points $[T_+]$ and $[T_-]$. Note that the stabilizers of $[T_+]$ and $[T_-]$ are virtually cyclic and virtually generated by $\varphi$ \cite{BestvinaFeignHandel1997}. Since $\varphi$ is atoroidal, $[T_+]$ and $[T_-]$ cannot be stabilized by any infinite order mapping class.

Projecting to the free splitting complex, we have that $[T_+]$ and $[T_-]$ are not finite distance apart. Since the free splitting complex is a coarse model for the sphere complex \cite{AramayonaSouto2011}, $[T_+]$ and $[T_-]$ are also not finite distance in the sphere complex. Consider a quasi-isometry of the sphere complex, for instance $\phi$,  that stabilizes these points. It cannot be bounded distance from any infinite order mapping class by the previous discussion. Moreover, $\phi$ cannot be bounded distance from a periodic mapping class because it is loxodromic. \end{proof}

Consider $(f_{\ast})^{-1}(G)$ in $\QI(\mathcal{A}(S_{g,1}))$. By \Cref{thm:rigidity}, $\QI(\mathcal{A}(S_{g,1})$ is itself isomorphic to $\Mod (S_{g,1})$. Since extended mapping class groups of finite-type surfaces are co-Hopfian \cite{BellMargalit2006, IvanovMcCarthy1999}, that is they are not isomorphic to any of their proper subgroups, it follows that $f_{\ast}^{-1}$ restricted to $G$ is surjective, but this precludes $f_{\ast}^{-1}$ from being an injection globally, a contradiction. \end{proof}

\section{QI-rigidity of the arc complex}
\label{sec:rigidity}

In this section we will prove \Cref{thm:rigidity}. We first remind the reader of some standard notation: given $\mathbb{N}$-valued functions $f, g$ over $\mathbb{N}$, by $f \prec g$ we will mean that there is some constant $D \ge 1$ so that 
\[ f(n) \leq D \cdot g(n) + D , \forall n \in \mathbb{N}.  \]
By $f \asymp g$ we will mean that $f \prec g$ and $g \prec f$.

Now, let $\phi: \mathcal{A}(S_{g,1}) \rightarrow \mathcal{A}(S_{g,1})$ be a quasi-isometry. Our goal will be to find some $f \in \mbox{Mod}^{\pm}(S_{g,1})$ and some $D>0$ so that for each arc $\alpha$, 
\[ d_{\mathcal{A}}(f(\alpha), \phi(\alpha)) < D. \]
To do this, we will first show that $\phi$ induces a quasi-isometry $\phi_{\ast}$ on the associated curve complex $\mathcal{C}(S)$. 

To construct $\phi_{\ast}$,  first select a complete and finite area hyperbolic metric on $S$. Let $S^{\circ} \subset S$ denote the complement of a standard cusp neighborhood about the single puncture, $p$; then there is some $K>0$ depending only on our choice of metric so that the diameter of $S^{\circ}$ is at most $K$. 

Now, fix some $x \in \mathcal{C}^{0}(S)$. Let $\gamma_{x}: [0,1] \rightarrow S^{\circ}$ denote a geodesic path of minimum length starting on the boundary of $S^{\circ}$ and ending at some point on the geodesic representative of $x$. Evidently, the hyperbolic length of $\gamma_{x}$ is at most $K$. Consider then the arc which we will denote by the concatenation $\gamma_{x} \ast x \ast \gamma_{x}^{-1}$ (where we have abused notation slightly and conflated $x$ with its geodesic representative): this is the arc that starts on $\partial S^{\circ}$, travels $\gamma_{x}$ to arrive at the geodesic representative of $x$, then travels along that representative, and finally travels back to $S^{\circ}$ along the inverse of $\gamma_{x}$. By extending this arc into the cusp neighborhood of $p$ by geodesic rays, we obtain a well-defined isotopy class of a proper simple arc on $S$, which we will denote by $a_{x}$. 

Consider then the arc $\phi(a_{x})$. Associated to it is a well-defined (up to isotopy) punctured annulus, obtained from $\phi(a_{x})$ by slightly pushing it first the left and then to the right, off of the puncture $p$ (these translates of $\phi(a_{x})$ will correspond to the two boundary components of the punctured annulus). Denote this punctured annulus by $A_{x}$. Finally, we will define $\phi_{\ast}(x)$ to be either component of $\partial A_{x}$.

\begin{remark} One could also define $\phi_{\ast}$ to be the full $2$-component multi-curve $\partial A_{x}$ and obtain a map into the power set of $\mathcal{C}^{0}$. Or, one could define $\phi_{\ast}$ in such a way that it sends an arc to the midpoint of the edge joining the curves in $\partial A_{x}$. Since we will only be concerned with coarse geometric invariants, these choices will not effect the argument. Having said that, the argument is slightly easier to state when we make the convention that $\phi_{\ast}$ is defined by simply choosing one of the two components of $\partial A_x$ at random and so this is what we will do going forward. 
\end{remark}

Because there may be several choices for the arc $\gamma_{x}$, we first need to argue that $\phi_{\ast}$ is coarsely well-defined. Indeed, there may be two distinct geodesic paths $\gamma, \gamma'$ that both minimize the distance from $\partial S^{\circ}$ to the geodesic representative of $x$. However, since the length of $\gamma$ and $\gamma'$ are both at most $K$, a standard collar lemma argument implies that they can only cross at most some $J = J(K)$ number of times. There is thus a universal upper bound on the geometric intersection number between the arcs $\gamma \ast x \ast \gamma^{-1}$ and $\gamma' \ast x \ast \gamma'^{-1}$. Therefore, the distance $d_{\mathcal{A}}$ in the arc graph between them is uniformly bounded by some constant $L$. 

Since $\phi$ is a quasi-isometry, the distance between the $\phi$-images of these arcs is at most some constant $L' = L'(L)$. Let $A, A'$ denote the punctured annuli corresponding to the $\phi$-images of the two arcs. Since the curve complex quasi-isometrically embeds in the arc and curve complex $\mathcal{AC}(S)$, we have 
\[ d_{\mathcal{C}(S)} (\partial A, \partial A') \prec d_{\mathcal{AC}(S)}( \partial A, \partial A') \leq 2 + L', \]
where the last inequality is just the triangle inequality together with the fact that $\partial A$ (resp. $\partial A'$) is disjoint from $\phi(\gamma \ast x \ast \gamma^{-1})$ (resp. $\phi(\gamma' \ast x \ast \gamma'^{-1})$). This proves coarse well-definedness. 

Next, we show that $\phi_{\ast}$ is coarsely Lipschitz: suppose $x, y$ are disjoint curves. Then by the same hyperbolic geometry argument used above, one obtains a universal upper bound (which we will again denote by $J = J(K)$) on the geometric intersection number between the corresponding arcs:
\[ i(\gamma_{x} \ast x \ast \gamma_{x}^{-1}, \gamma_{y} \ast y \ast \gamma_{y}^{-1}) < J. \]
Hence, their distance in the arc graph is uniformly bounded and so their $\phi$-images are uniformly close as well. The argument now concludes in an identical fashion to the last paragraph.

Finally, we will show that $\phi_{\ast}$ admits a coarsely Lipschitz coarse inverse map, which in particular implies that $\phi_{\ast}$ is coarsely surjective and that it contracts distances by only some universally bounded amount, completing the proof that it is a quasi-isometry. Given $x \in \mathcal{C}^{0}$, define a new curve $\psi(x)$ by the same procedure as with $\phi_{\ast}$, only now we replace $\phi$ with $\phi^{-1}$: convert $x$ into an arc by concatenation with short paths to $\partial S^{\circ}$, apply $\phi^{-1}$, and then consider the corresponding punctured annulus-- choose one of its two boundary components at random. 

By the exact same arguments used above $\psi$ is coarsely Lipschitz. So it suffices to prove that there is some universal constant $B$ so that 
\[ d_{\mathcal{A}}(x, \psi \phi_{\ast}x) \leq B.\]

The key observation is that, while we do not have control over the isotopy class of the arc $\gamma_{A_{x}}$ rel $\partial S^{\circ} \cup A_{x}$ (where $\gamma_{A_x}$ is formed using either component of $\partial A_{x}$ in the same way as $\gamma_x$ was previously), we do know that this arc will be contained in $A_{x}$ simply because $ \partial A_{x}$ separates $\partial S^{\circ}$ from the rest of the surface. Therefore, one easily sees that the arc $\gamma_{A_{x}} \ast \phi_{\ast}x \ast \gamma_{A_{x}}^{-1}$ is homotopic to the arc $\phi(x)$. Therefore, 
\[ d_{\mathcal{A}} (\phi^{-1} (\phi_{\ast} (x), \gamma_{x} \ast x \ast \gamma_{x}^{-1}) < B' \]
for some universal constant $B'$, whence it follows that 
\[ d_{\mathcal{C}} (\psi \phi_{\ast}(x) , x) \leq B' + 2 =: B, \]
as desired.

Now, by Rafi--Schleimer \cite{RafiSchleimber2011}, the map $\phi_{\ast}$ is bounded distance from some mapping class $f$: there is $C >0$ so that for all $x \in \mathcal{C}^{0}$, 
\[ d_{\mathcal{C}} (\phi_{\ast}(x) , f(x)) < C. \]
Our goal is to show the analogous statement for arcs: 
\[ d_{\mathcal{A}}(\phi(\alpha), f(\alpha)) = O(1), \forall \alpha \in \mathcal{A}^{0}. \]

For this, we will rely on the familiar distance formula for the arc graph (see \cite{MasurSchleimer2013} for details):
\[ d_{\mathcal{A}}(\alpha, \beta) \asymp \sum_{Y \in \Omega(S)} [[d_{Y}(\alpha, \beta]]_{R}, \]
where: 
\begin{itemize}
\item $\Omega(S)$ denotes the set of \textit{witness} subsurfaces (which in this case is any subsurface containing the puncture $p$);
\item $d_{Y}(\alpha, \beta)$ denotes the distance in $\mathcal{C}(Y)$ between the subsurface projections of $\alpha, \beta$ to $Y$; 
\item $R>0$ is a constant depending only on $\chi(S)$; and 
\item $[[ x ]]_{R}= x$ whenever $x \geq R$ and is $0$ otherwise.
\end{itemize}

It therefore suffices to prove that there is some universal $T>0$ so that  for $Y \in \Omega(S)$, 
\[ d_{Y}(\phi(\alpha), f(\alpha)) < T. \]
Note that  $d_{Y}(f(\alpha), f(\partial A_{\alpha}))$ is uniformly bounded-- let us say by some constant $T_{1}$-- because $\alpha$ and $\partial A_{\alpha}$ are disjoint, $f$ is a homeomomorphism and $\pi_{Y}$ is coarsely Lipschitz (see \cite{MasurMinsky00} for details). Moreover, $d_{Y}(\phi(\alpha), \phi_{\ast}(\partial A_{\alpha}))< T_{2}$ for some other universal constant $T_{2}$, because $\phi_{\ast}$ is defined in a way which guarantees that $\phi_{\ast}(\partial A_{\alpha})$ will be disjoint from $\phi(\alpha)$. 

Now, 
\begin{align*} 
 d_{Y}(&\phi(\alpha), f(\alpha)) \\ 
 & \leq d_{Y}(\phi(\alpha), \phi_{\ast}(\partial A_{\alpha}))+ d_{Y}(\phi_{\ast}(\partial A_{\alpha}), f(\partial A_{\alpha}))+ d_{Y}(f(\partial A_{\alpha}),f(\alpha)) \\ 
& < T_{1} + C + T_{2} = : T.
\end{align*}

This completes the proof. 

\section*{Acknowledgments} The first author thanks Sam Taylor for helpful conversations. The second author thanks Caglar Uyanik for his invaluable expertise on the dynamics of $\Out(F_{n})$; in particular, for sketching the proof of the claim used in the proof of \Cref{thm:main}. Both authors acknowledge support from the National Science Foundation through grants DMS-1939936 (Aougab) and DMS-2231286 (Loving).

\bibliographystyle{alpha}
\bibliography{main} 

\end{document}